\documentclass[11pt]{article}
\usepackage{amssymb}
\usepackage{amsmath}

\newtheorem{thm}{Theorem}[section]
\newtheorem{lemma}{Lemma}[section]

\newtheorem{remark}{Remark}[section]
\newcommand{\R}{\mathbb{R}}

\newcommand{\N}{\mathbb{N}}

\def\R{\mathbb R}

\def\C{\complement}
\newcommand{\eps}{\varepsilon}
\newcommand{\vphi}{{\varphi}}
\def\ds{\displaystyle}
\def\N{\mathbb N}
\def\div{\mbox{div}\, }

\def\na{\nabla}
\def\pa{\partial}

\title{On the isentropic compressible Navier-Stokes equation}

\author{A.Mellet\thanks{mellet@math.utexas.edu}
  and A.Vasseur\thanks{vasseur@math.utexas.edu}\\
 \\
\small Department of Mathematics\\
\small University of Texas at Austin
 }

\date{}

\begin{document}
\maketitle
\begin{abstract}
In this article, we consider  the compressible  Navier-Stokes
equation with density dependent viscosity coefficients. We focus on
the case where those coefficients vanish on vacuum. We prove the
stability of weak solutions   for  periodic domain $\Omega = T^{N}$
as well as the whole space $\Omega = \R^{N}$, when $N=2$ and $N=3$.
The pressure is given by $p=\rho^{\gamma}$, and our result holds
for any $\gamma>1$. In particular, we prove the stability of weak
solutions of the Saint-Venant model for shallow water.
\end{abstract}

\section{Introduction}

This paper is devoted to the Cauchy problem of the compressible
Navier-Stokes equation with viscosity coefficients vanishing on
vacuum. Let $\rho(t,x)$ and $u(t,x)$ denote the density and the
velocity of an isentropic compressible viscous fluid (as usual,
$\rho$ is a non-negative function and $u$ is a vector valued
function, both defined on a subset $\Omega$ of $\R^{N}$). Then, the
Navier-Stokes equation for isentropic  compressible viscous fluids
reads (see \cite{LL}):
\begin{equation}\label{eq:0}
\begin{array}{l}
 \pa_{t} \rho +\div (\rho u) = 0 \\
 \pa_{t}(\rho u) +\div (\rho u \otimes u) + \nabla_{x} p -\div (h\, D(u)) - \na (g\, \div u)= 0
\end{array}
\end{equation}
where $p(\rho)=\rho^\gamma$, $\gamma>1$, denotes the pressure,
$D(u)=\frac{1}{2}[\na u+{}^{t}\na u]$ is the strain tensor and $h$
and $g$ are the  two Lam\'e viscosity coefficients (depending on the
density $\rho$) satisfying
\begin{equation}\label{eq:h+Ng}
 h>0 \qquad h+Ng \geq 0
\end{equation}
($h$ is sometime called the shear viscosity of the fluid, while $g$
is usually referred to as the second viscosity coefficient).
One of the major difficulty of compressible fluid mechanics is to
deal with vacuum. The problem of existence of global solution in
time for Navier-Stokes equations was addressed in one dimension for
smooth enough data by Kazhikov and Shelukhin \cite{K}, and for
discontinuous one, but still with densities away from zero, by Serre
\cite{Serre} and Hoff \cite{Hoff}. Those results have been
generalized to higher dimensions by Matsumura and Nishida
\cite{nishida} for smooth data close to equilibrium and by Hoff
\cite{Hoff2}, \cite{Hoff3} in the case of discontinuous data.

Concerning large initial data, Lions showed in \cite{PLL2} the
global existence of weak solutions for $\gamma\geq3/2$ for $N=2$ and
$\gamma\geq 9/5$ for $N=3$. This result has been extended later by
Feireisl, Novotny, and Petzeltova   to the range $\gamma>3/2$ in
\cite{F2}, and very recently by Feireisl to the full system of the
Navier-Stokes equations involving the energy equation \cite{F3}.
Other results provide the full range $\gamma>1$ under symmetries
assumptions on the initial datum (see for instance Jiang and Zhang
\cite{J1}). 
All those results do not require to be far from the vacuum.
However they rely strongly on the assumption that the viscosity
coefficients are bounded below by a positive constant. This non
physical assumption allows to get some estimates on the gradient of
the velocity field.

The main difficulty when dealing with vanishing viscosity
coefficients on vacuum is that the velocity cannot even be defined
when the density vanishes.
 The first result handling this difficulty is due to
 Bresch, Desjardins and Lin \cite{BDL}.  They showed the
$L^1$ stability of weak solutions for the following Korteweg's system of
equations:
\begin{equation}\label{eq:korteweg}
\begin{array}{l}
\pa_{t} \rho +\div (\rho u) = 0 \\
 \pa_{t}(\rho u) +\div (\rho u \otimes u) + \nabla_{x} p -\nu \div (\rho \, D(u)) = \kappa \rho \na \Delta \rho.
 \end{array}
\end{equation}
The result was later improved by Bresch and Desjardins in \cite{BD2} to include the case of vanishing capillarity ($\kappa=0$), but with an additional quadratic friction term $r \rho |u|u$
 (see also \cite{BD}). 
The key point in those papers is to show that the structure of the
diffusion term provides some regularity for the density thanks to a new entropy inequality.
However, those estimates are not enough to treat 
the case without capillarity and friction effects $\kappa=0$ and $r=0$ (which
corresponds to equation (\ref{eq:0}) with $h(\rho)=\rho$ and
$g(\rho)=0$). 

The main difficulty, to prove the stability of the solutions of (\ref{eq:0}), 
is to pass to the limit in the term $\rho(u\otimes u)$ (which requires the strong convergence of $\sqrt{\rho}u$).
Note that this is easy when the viscosity coefficients are bounded below by a positive constant.
On the other hand, the new bounds on the gradient of the density make the control of the pressure term far simpler than in the case of constant viscosity coefficients.

Our result is in the same spirit as the one of Bresch, Desjardins and Lin and makes use of the same entropy inequality, first discovered by  Bresch and Desjardins in \cite{BD} for the particular case where $h(\rho)=\rho$ and $g(\rho)=0$.
We actually  use a slightly more general estimate, which holds for any viscosity coefficients $h(\rho),g(\rho)$ satisfying the relation:
\begin{equation}
g(\rho) = \rho h'(\rho) - h(\rho).  \label{ghanticipee}
\end{equation} 
This estimate first appeared in a Note 
by Bresch and Desjardins \cite{BD3} in the context of Korteweg systems of equations.
However, we will see that the capillary term is by no means necessary to the derivation of the crucial estimates which thus hold for the compressible Navier-Stokes system (\ref{eq:0}).

Our main contribution is to show the $L^1$
stability of weak solutions of (\ref{eq:0}) 
under some conditions  on the viscosity coefficients (including
(\ref{ghanticipee})) but without any additional regularizing terms.
The interest of our result lie primarily in the fact that our conditions allow for viscosity coefficients that vanish on the vacuum set.
It includes the case $h(\rho)=\rho$, $g(\rho)=0$ 
(when $N=2$ and $\gamma=2$, we recover the  Saint Venant model for
Shallow water), 
but our conditions on $h$ and $g$ will exclude the case of 
constant viscosity $h(\rho)=\mu$, $g(\rho)=\xi$. Indeed, it is readily seen that
 (\ref{ghanticipee}) implies that $g(\rho)=\xi=-\mu$, and
thus $\mu+\xi=0$. 
In this border line case we thus lose all informations on the derivatives of $u$.
It is worth pointing out that while we can gain regularity on the density with this new estimate, we have to loose regularity on the velocity (on the vacuum set).

Note that the main difficulty will be to establish the compactness of $\sqrt{\rho}u$ in $L^{2}$ strong, and
the key ingredient to achieve this is an additional estimate which  bounds
$\sqrt{\rho}u$ in $L^{\infty}(0,T;L^{2+2\alpha}(\Omega))$ for some small $\alpha>0$
(the usual entropy estimate only gives a bound in $L^{\infty}(0,T;L^{2}(\Omega))$).

For the sake of simplicity we will consider the case $\Omega=\R^N$
and the case of bounded domain with periodic boundary conditions,
namely $\Omega=T^N$. For the same reason we consider only power
pressure laws although the result could be extend to non monotonic
pressure law of the form of \cite{F1}. 
Note that the result holds
for any power $\gamma>1$ under appropriate assumptions on $h$ and
$g$. Classically, $L^1$ stability  is considered as the main step to
prove the existence of weak solutions. 
To obtain the existence of weak solutions, one is thus left with the
technical task of constructing a sequence of approximated solutions
verifying the a priori estimates. 
Although this final step is in
most cases quite standard, we point out that in this
particular situation it seems highly non trivial because of the complexity of the
additional entropy inequality.

\vspace{5pt}

In the next section, we state the  assumptions on the viscosity
coefficients, define precisely the notion of ``weak solutions'' and
state our main results. In Section \ref{sec:est}, we recall the well
known physical energy inequality and state the key  estimates. The proof
of Theorem \ref{thm:stability} is detailed in Section
\ref{sec:proof}.
For the sake of completeness, we give in Section \ref{sec:estimate} the proof of the  entropy inequality of Bresch and Desjardins in the context of compressible Navier-Stokes equation.

\vspace{10pt}

\section{Notations and main result}
Let $\Omega$ denote a subset of $\R^{N}$.
We assume that  $\Omega $ is either the whole space $\R^{N}$ or a bounded domain with periodic boundary conditions ($\Omega = T^{N}$).
For the sake of simplicity, we will take $D(u) =\na u$, though the full strain tensor could be considered  without any additional difficulty.
This leads to the following system of equations:
\begin{eqnarray}
&& \pa_{t} \rho +\div (\rho u) = 0 \label{eq:euler1}\\
&& \pa_{t}(\rho u) +\div (\rho u \otimes u) + \nabla_{x} \rho^{\gamma} -\div (h(\rho)\na u) - \na (g(\rho)\div u)= 0, \label{eq:euler2}
\end{eqnarray}
with initial conditions
\begin{equation}\label{eq:initial3}
\rho|_{t=0} =\rho_{o} \geq 0\, ,\qquad \rho u|_{t=0} = m_{o}.
\end{equation}
Before introducing the notion of weak solution, let us state the
assumptions we make on the viscosity coefficients.

\vspace{20pt}

\noindent {\bf Conditions on $h(\rho)$ and $g(\rho)$:} \\
First we consider $f(\rho),g(\rho)$  verifying:
\begin{equation}\label{gh}
g(\rho) = \rho h'(\rho)-h(\rho).
\end{equation}
As stated in the introduction, this structure constraint is
fundamental to get more regularity on the density. Moreover, we
assume that there exists a positive constant $\nu\in(0,1)$ such that
\begin{eqnarray}
&& h'(\rho)\geq \nu\, ,  \qquad h(0)\geq 0\label{h}\\
&& |g'(\rho)| \leq \frac{1}{\nu} h'(\rho) \label{g}\\
&& \nu h(\rho)\leq h(\rho)+N g(\rho) \leq \frac{1}{\nu} h(\rho) \,
.\label{g+h}
\end{eqnarray}
When $\gamma \geq 3$ and $N=3$, we also require that
\begin{equation} \label{hg}
\liminf_{\rho\rightarrow \infty}
\frac{h(\rho)}{\rho^{\gamma/3+\eps}} >0,
\end{equation}
for some small $\eps>0$.
\vspace{5pt} Let us make some remarks about those assumptions.
\begin{remark}
The functions
$$ h(\rho) = \rho,\qquad g(\rho)=0$$
satisfy (\ref{gh}-\ref{g+h}). In fact, any linear combination of
$\rho^{k}$ with $k\geq 1$ is an admissible function for $h(\rho)$.
\end{remark}

\begin{remark}
The lower estimate in (\ref{g+h}) is trivial when $g\geq 0$, while
the upper estimate is trivial  when $g\leq 0$. Together they yield:
$$ |g(\rho)| \leq C_{\nu} h(\rho)\, \quad \forall \rho>0.$$
This inequality and (\ref{g}) will be necessary to pass to the limit
in the term $\na (g(\rho_{n})\div u_{n})$.
\end{remark}

\begin{remark} \label{rmk:1}
Condition (\ref{h}) makes the proof simpler, but is not optimal.
However, condition (\ref{g+h}) is necessary to control the viscosity
term and together with (\ref{gh}), it yields
$$\frac{N-1+\nu}{N\rho} \leq \frac{h'(\rho)}{h(\rho)}\leq \frac{N-1+1/\nu}{N\rho}, \quad \mbox{ for all $\rho> 0$,}$$
and so
\begin{equation} \label{eq:h2/3}
\left\{
\begin{array}{ll}
C \rho ^{(N-1)/N+\nu/N} \leq  h(\rho) \leq C \rho ^{(N-1)/N+1/(N\nu)}  , \quad & \rho \geq 1\\
C \rho ^{(N-1)/N+1/(N\nu)} \leq h(\rho) \leq C \rho^{(N-1)/N+\nu/N},
& \rho \leq 1
\end{array}
\right.
\end{equation}
In particular, we must have $h(0)=0$. Moreover, this shows that if
we do not assume (\ref{h}), the ``best'' $h(\rho)$ we can take is
$h(\rho) = \rho^{(N-1)/N+\nu/N}$. This is actually enough to prove
the stability  of weak solutions for all $\gamma$ when $N=2$ and
for $\gamma < 3/2$ when $N=3$. However, if we assume $h(\rho) \sim
C\rho^{2/3+\nu}$ for small $\rho$ and $h(\rho) \sim C \rho$ for
large $\rho$, then we can take any $\gamma\in(1,3)$ when $N=3$.
\end{remark}

\vspace{30pt} 

\noindent{\bf Notion of weak solutions}

 We say that $(\rho,u)$ is a weak solution of
(\ref{eq:euler1}-\ref{eq:euler2}) on $\Omega\times [0,T]$, with
initial conditions (\ref{eq:initial3}) if
\begin{equation*} 
\begin{array}{l}
\ds \rho\in L^{\infty}(0,T, L^{1}(\Omega)\cap L^{\gamma}(\Omega)) , \\
\sqrt \rho \in L^{\infty}(0,T;H^{1}(\Omega)) ,\\
\ds \sqrt \rho\, u \in L^{\infty} (0,T; (L^{2}(\Omega))^{N}) , \\
\ds h(\rho) \na u  \in
L^{2}(0,T;(W^{-1,1}_{\mathrm{loc}}(\Omega))^{N\times N}), \quad
g(\rho) \div u \in L^{2}(0,T;W^{-1,1}_{\mathrm{loc}}(\Omega)),
\end{array}
\end{equation*}
with $\rho\geq 0$ and $(\rho,\sqrt{\rho} u)$ satisfying
$$
\left\{
\begin{array}{l}
 \pa_{t} \rho+ \div (\sqrt \rho \sqrt \rho u) = 0 \\ \rho(0,x)= \rho_{o}(x) \end{array}
\right. \qquad \mbox{ in } \mathcal D',$$
and if the following equality holds for  all $\vphi(t,x)$ smooth test function with compact support such that $\vphi(T,\cdot)=0$:
\begin{eqnarray}
&&\!\!\!\!\!\!\!\!\!\!\!\! \int_{\Omega} m_{o} \cdot \vphi(0,\cdot)\, dx +
\int_{0}^{T}\!\! \!\!\int_{\Omega} \sqrt \rho (\sqrt \rho u ) \pa_t \vphi +\sqrt \rho u \otimes \sqrt \rho u : \nabla \vphi\, dx \nonumber\\
&&\!\!\!\!\!\!\!\!\!\!\!\! \quad + \int_{0}^{T}\!\!\!\!\int_{\Omega} \rho^{\gamma}\div \vphi \, dx   - \big\langle h(\rho) \na u \, ,\,  \na \vphi\big\rangle  - \big\langle g(\rho) (\div u)\, ,\, (\div \vphi)\big\rangle  =0 ,\label{eq:weak}
\end{eqnarray}
where the diffusion terms make sense when written as
\begin{eqnarray*}
&&  \big\langle h(\rho) \na u\, ,\,  \na \vphi\big\rangle =\\
&& \qquad = -\int \frac{h(\rho)}{\sqrt\rho}(\sqrt \rho u_{j}) \pa_{ii} \vphi_{j} \, dx\, dt - \int (\sqrt\rho u_{j}) 2h'(\rho) \pa_{i}\sqrt\rho \pa_{i} \vphi_{j} \, dx \, dt,
\end{eqnarray*}
and
\begin{eqnarray*}
&&  \big\langle g(\rho) (\div u) \,,\,  (\div \vphi)\big\rangle= \\
&& \qquad = -\int \frac{g(\rho)}{\sqrt\rho}(\sqrt \rho u_{i}) \pa_{ij} \vphi_{j} \, dx\, dt - \int (\sqrt\rho u_{i}) 2g'(\rho) \pa_{i}\sqrt\rho \pa_{j} \vphi_{j} \, dx \, dt.
\end{eqnarray*}
In particular, the fact that the diffusion term $h(\rho) \na u$ (and
$g(\rho) \div u$) lies in $
L^{2}(0,T;(W^{-1,1}_{\mathrm{loc}}(\Omega))^{n\times n})$ will
follow from the fact that
$$ h'(\rho) \na\sqrt\rho \in L^{\infty}(0,T;L_{\mathrm{loc}}^{2}(\Omega))\, ,
\quad\mbox{and} \quad h(\rho)/\sqrt\rho \in L^{\infty}(0,T;L^{2}_{\mathrm{loc}}(\Omega)) ,$$
and similar conditions on $g(\rho)$. This will be provided by
assumptions (\ref{g}), (\ref{h}) and (\ref{eq:h2/3}).

\vspace{20pt} \noindent{\bf Main result:}

 The main result of this
paper is the following:
\begin{thm}\label{thm:stability}
Assume that $\gamma>1$ and that $h(\rho)$ and $g(\rho)$ are two
$C^{2}$ functions of $\rho$ satisfying  conditions
(\ref{gh})-(\ref{g+h}) (together with (\ref{hg}) if $\gamma \geq 3$
and $N=3$). Let $(\rho_{n},u_{n})_{n\in \N}$ be a sequence of weak
solutions of (\ref{eq:euler1}-\ref{eq:euler2}) satisfying entropy
inequalities (\ref{eq:entropy1}), (\ref{eq:entropy2}) and
(\ref{eq:moment}), with initial data
$$\rho_{n}|_{t=0} = \rho^{n}_{o}(x)\quad \mbox{ and }\quad \rho_{n}u_{n}|_{t=0} = m^{n}_{o}(x)=\rho^{n}_{o}(x)u^{n}_{o}(x),$$
where $\rho^{n}_{o}$ and $u^{n}_o$ are such that
\begin{equation}\label{initial00}
\rho^{n}_{o}\geq 0,\qquad \rho^{n}_{o} \rightarrow \rho_{o}\mbox{ in } L^{1}(\Omega) , \quad \rho_{o}^{n}u^{n}_{o}\rightarrow \rho_{o}u_{o} \mbox{ in } L^{1}(\Omega)  ,
\end{equation}
and satisfy the following bounds (with $C$ constant independent on $n$):
\begin{equation}\label{initial11}
\int_{\Omega} \rho^{n}_{o} \frac{|u^{n}_{o}|^{2}}{2} + \frac{1}{\gamma-1}{\rho^{n}_{o}}^{\gamma}\, dx < C,
\qquad
\int_{\Omega}\frac{1}{ \rho^{n}_{o}}|\na h(\rho^{n}_{o})|^{2}\, dx <C,
\end{equation}
and
\begin{equation}\label{initial22}
\int_{\Omega} \rho^{n}_{o} \frac{|u^{n}_{o}|^{2+\delta}}{2} \, dx < C,
\end{equation}
for some small $\delta>0$.

Then, up to a subsequence, $(\rho_{n},\sqrt{\rho_{n}}u_{n})$
converges strongly to a weak solution of
(\ref{eq:euler1})-(\ref{eq:euler2}) satisfying entropy inequalities
(\ref{eq:entropy1}), (\ref{eq:entropy2}) and (\ref{eq:moment}) (the
density $\rho_{n}$ converges strongly in $\mathcal
C^{0}((0,T);L_{loc}^{3/2}(\Omega))$, $\sqrt{\rho_{n}}u_{n}$
converges strongly in $L^2(0,T;L^2_{\mathrm{loc}}(\Omega))$ and the
momentum $m_{n}=\rho_{n}u_{n}$ converges strongly in $L^{1}(0,T;
L^1_{\mathrm{loc}}(\Omega))$, for any $T>0$).
\end{thm}

\section{Entropy inequalities and a priori estimates}\label{sec:est}
In this section, we recall the well-known energy inequality and
state the main inequalities that we will use throughout the proof of
Theorem \ref{thm:stability}.

The usual energy inequality associated with the system of equations (\ref{eq:euler1}-\ref{eq:euler2}) can be written as:
\begin{equation}\label{eq:entropy1}
 \frac{d}{dt}\int \rho \frac{u^{2}}{2} + \frac{1}{\gamma-1} \rho^{\gamma} \,dx +
\int h(\rho)|\na u|^{2}\, dx  + \int g(\rho) (\div u)^{2}\, dx \leq  0.
\end{equation}
This inequality can be established for smooth solutions of
(\ref{eq:euler1}-\ref{eq:euler2})   by multiplying the momentum
equation by $u$.

When  $h$ and $g$ satisfies $h(\rho) + N g(\rho)\geq 0$ and if the initial data are taken in such a way that
$$
\mathcal E_{o} = \int_{\Omega} \rho_{o} \frac{u_{o}^{2}}{2} + \frac{1}{\gamma-1}\rho_{o}^{\gamma}\, dx < +\infty,
$$
then (\ref{eq:entropy1}) yields:
\begin{equation} \label{eq:bounds1}
\begin{array}{l}
\ds || \sqrt \rho\, u ||_{L^{\infty}(0,T;L^{2}(\Omega))} \leq C, \\
\ds || \rho||_{L^{\infty}(0,T;L^{\gamma}(\Omega))} \leq C .
\end{array}
\end{equation}
Furthermore, Hypothesis (\ref{g+h}) gives:
\begin{equation} \label{eq:bounds1h}
 ||\sqrt{h(\rho)} \na u||_{L^{2}(0,T;L^{2}(\Omega))} \leq C.
\end{equation}
Finally, integrating (\ref{eq:euler1}) with respect to $x$ yields the natural $L^{1}$ estimate:
$$ ||\rho ||_{L^{\infty}(0,T;L^{1}(\Omega))}\leq C.$$

Unfortunately, it is a well-known fact that those estimates are not enough to prove the stability of the solutions of  (\ref{eq:euler1}-\ref{eq:euler2}).
In particular, the fact that $\rho^{\gamma}$ is bounded in $L^{\infty}(0,T;L^{1}(\Omega))$ does not implies that  $\rho_{n}^{\gamma}$ converges to $\rho^{\gamma}$.
\vspace{15pt}

However, further estimates can be obtained by mean of the following lemma (the proof of which is postponed to Section \ref{sec:estimate}):
\begin{lemma}\label{lem:estimate}
Assume that $h(\rho)$ and $g(\rho)$ are two $C^{2}$ functions such
that (\ref{gh}) holds true.
Then, the following inequality holds for  smooth solutions of (\ref{eq:euler1}-\ref{eq:euler2}):
\begin{equation}\label{eq:entropy2}
\frac{d}{dt}\int \frac{1}{2}\rho |u+\na \vphi(\rho)|^{2}+   \frac{1}{\gamma-1} \rho^{\gamma} \,dx +
\int \na \vphi(\rho) \cdot \nabla \rho^{\gamma}\, dx \leq 0,
\end{equation}
with $\vphi$ such that
\begin{equation}\label{phi}
 \vphi' = \frac{h'}{\rho}.
 \end{equation}
\end{lemma}
This lemma is similar to the result of D.~Bresch and B. Desjardin in \cite{BD3}, in which the same inequality was derived when capillary effects are taken into account.

We immediately see that since the viscosity coefficient $h(\rho)$ is
an increasing function of $\rho$ and when the initial data satisfies
$$
\int_{\Omega} \rho_{o}|\na \vphi(\rho_{o})|^{2}\, dx <+\infty,
$$
inequality (\ref{eq:entropy2}) yields:
\begin{equation}\label{eq:bounds2}
\frac{1}{2}|| \sqrt \rho \na \vphi(\rho)||_{L^{\infty}(0,T;L^{2}(\Omega))} = || h'(\rho) \na \sqrt \rho||_{L^{\infty}(0,T;L^{2}(\Omega))}\leq C,
\end{equation}
and
\begin{equation}\label{eq:bounds2ter}
||\sqrt{h'(\rho)\rho^{\gamma-2}}\na \rho||_{L^{2}(0,T;L^{2}(\Omega))} \leq C.
\end{equation}
Under assumption (\ref{h}) on  $h$, those estimates will give
additional control on the density $\rho$ and on the pressure
$\rho^\gamma$, which will be enough to prove the stability of weak
solution. \vspace{10pt}

\vspace{15pt}

Finally, we shall make use of the following result:

\begin{lemma}\label{lem:estimate2}
Assume
$$ h(\rho) + N g(\rho) \geq \nu h(\rho)$$
for some $\nu \in (0,1)$ (which is a part of (\ref{g+h})), and let
$\delta\in(0,\nu/4)$. Then, smooth solutions of
(\ref{eq:euler1}-\ref{eq:euler2}) satisfy the following inequality:
\begin{equation}\label{eq:moment}
\begin{array}{l}
\ds{\frac{d}{dt} \int \rho \frac{|u|^{2+\delta}}{2+\delta}\,dx +
\frac{\nu}{4} \int h(\rho)|u|^{\delta} |\na u|^{2} \, dx}\\[0.5cm]
\qquad \qquad\ds{\leq
 \left( \int \left(\frac{\rho^{2\gamma-\delta/2}}{h(\rho)}\right)^{2/(2-\delta)}\, dx\right)^{(2-\delta)/2} \left(\int \rho |u|^{2}\,dx
 \right)^{\delta/2}.}
\end{array}
\end{equation}
where $|\na u|^{2} = \sum_{i}\sum_{j}|\pa_{i } u_{j}|^{2}$.
\end{lemma}
This inequality is quite simple to establish and will be essential
in the proof of Theorem \ref{thm:stability} to prove that
$\sqrt{\rho_{n}}u_{n}$ is bounded in
$L^{\infty}(0,T;L^{2+2\alpha}(\Omega))$ (see Lemma
\ref{lem:moment}). Note, however, that to derive further estimates
from this inequality, we need to control the right hand side of
(\ref{eq:moment}). Inequality (\ref{eq:entropy1}) immediately
provide a bound on $\int \rho |u|^{2}\,dx$, so the problem will be
to control enough power of $\rho$ to get a bound on $\int
\left(\frac{\rho^{2\gamma-\delta/2}}{h(\rho)}\right)^{2/(2-\delta)}\,
dx$. This will be achieved using  (\ref{eq:bounds2ter}). Of course,
we also need to assume that the initial condition satisfies
$$
 \int \rho_{o} \frac{ |u_{o}|^{2+\delta}}{2} \, dx < C .
$$

{\it Proof of Lemma \ref{lem:estimate2}.}
Let $\delta\in (0,\nu/4)$.
Multiplying (\ref{eq:euler2}) by $u|u|^{\delta}$, we get:
\begin{eqnarray*}
&&\int  \rho \pa_{t} \frac{|u|^{2+\delta}}{2+\delta}\,dx + \int \rho u\cdot \na \frac{|u|^{2+\delta}}{2+\delta}\,dx \\
&& \qquad +\int h(\rho)|u|^{\delta} (\na u)^{2}\,dx + \delta \int h(\rho) |u|^{\delta-2} u_{i} u_{k} \pa_{j} u_{i} \pa_{j} u_{k}\,dx \\
&& \qquad +\int g(\rho) |u|^{\delta} (\div u)^{2}\,dx + \delta \int g(\rho) |u|^{\delta -2} u_{k}u_{j}\pa_{i} u_{i}\pa_{j} u_{k}\,dx  \\
&& \qquad + \int |u|^{\delta } u \cdot \na \rho^{\gamma}\, dx =0.
\end{eqnarray*}
Since
$$
 (\div u)^{2} = \sum_{i} \sum_{j} \pa_{i} u_{i}\pa_{j}u_{j} \leq \sum_{i} \sum_{j} \frac{1}{2}(\pa_{i} u_{i}^{2}+\pa_{j}u_{j}^{2})  \leq N |\na u |^{2},
$$
condition (\ref{g+h}) yields:
\begin{eqnarray*}
&&\int  \rho \pa_{t} \frac{|u|^{2+\delta}}{2+\delta}\,dx + \int \rho
u\cdot \na \frac{|u|^{2+\delta}}{2+\delta}\,dx  +\nu \int
h(\rho)|u|^{\delta} (\na u)^{2}\,dx\\
&&\qquad\qquad\qquad
+ \int |u|^{\delta } u \cdot \na \rho^{\gamma}\, dx \\
&& \qquad \qquad \leq  \delta \int h(\rho) |u|^{\delta-2} u_{i}
u_{k} \pa_{j} u_{i} \pa_{j} u_{k}\,dx \\
&&\qquad\qquad\qquad+ \delta \int g(\rho) |u|^{\delta -2}
u_{k}u_{j}\pa_{i} u_{i}\pa_{j} u_{k}\,dx ,
\end{eqnarray*}
and since $\delta<\nu/4$,  we deduce:
\begin{eqnarray*}
&& \int \rho \pa_{t}  \frac{|u|^{2+\delta}}{2+\delta}\,dx + \int \rho u\cdot \na \frac{|u|^{2+\delta}}{2+\delta}\,dx  +\frac{\nu}{2} \int h(\rho)|u|^{\delta} (\na u)^{2}\,dx \\
&& \qquad\qquad\qquad\qquad\qquad\qquad\qquad\qquad\qquad  + \int |u|^{\delta } u \cdot \na \rho^{\gamma}\, dx \leq 0.
\end{eqnarray*}
Moreover, multiplying (\ref{eq:euler1}) by $ \frac{|u|^{2+\delta}}{2+\delta}$ and integrating by parts, we have
$$
\int  \frac{|u|^{2+\delta}}{2+\delta} \pa_{t } \rho\,dx - \int \rho
u \cdot \na  \frac{|u|^{2+\delta}}{2+\delta}  \, dx =0
$$
and summing  the last two inequalities, we get:
$$
\frac{d}{dt} \int \rho \frac{|u|^{2+\delta}}{2+\delta}\,dx +
\frac{\nu}{2} \int h(\rho)|u|^{\delta} |\na u|^{2} \, dx \leq \left|
\int |u|^{\delta} u\cdot \na \rho^{\gamma}\, dx \right|,
$$

It remains to bound  the right hand side. We have:
\begin{eqnarray*}
&&\left| \int |u|^{\delta} u\cdot \na \rho^{\gamma}\, dx \right|  =  \left|  -\int \rho^{\gamma} |u|^{\delta} \div u \,dx- \delta \int \rho^{\gamma} |u|^{\delta-2} u(u\cdot\na) u   \,dx \right|\\
&&\qquad\qquad \leq  (\sqrt N+\delta) \left| \int \rho^{\gamma} |u|^{\delta} |\na u| \,dx \right|\\
&&\qquad\qquad \leq  (\sqrt N+\delta) \left(\int h(\rho) |u|^{\delta} |\na u|^{2}\, dx \right)^{1/2} \left(\int \frac{\rho^{2\gamma}}{h(\rho)} |u|^{\delta} \,dx \right)^{1/2}\\
&&\qquad\qquad \leq  \frac{\nu}{4} \int h(\rho) |u|^{\delta} |\na
u|^{2}\, dx + C_{\nu} \int \frac{\rho^{2\gamma}}{h(\rho)}
|u|^{\delta}\, dx,
\end{eqnarray*}
where  the last term satisfies (if $\delta\in (0,2)$):
$$\int \frac{\rho^{2\gamma}}{h(\rho)} |u|^{\delta}\, dx \leq \left( \int \left(\frac{\rho^{2\gamma-\delta/2}}{h(\rho)}\right)^{2/(2-\delta)}\, dx\right)^{(2-\delta)/2} \left(\int \rho |u|^{2}\,dx \right)^{\delta/2},
$$
and the lemma follows.

\vspace{20pt}

We now have all the necessary tools to  prove Theorem \ref{thm:stability}.

\section{Proof of Theorem \ref{thm:stability}} \label{sec:proof}
We now present the proof of Theorem \ref{thm:stability}. To begin
with, we need to make precise the assumptions on the initial data.
\vspace{10pt}

\noindent {\bf Initial data:} \\
We recall that the initial data must satisfy (\ref{initial11}), and
(\ref{initial22}) to make use of all the inequalities presented in
the previous section:
\begin{equation} \label{eq:initial}
\begin{array}{l}
\displaystyle\rho^{n}_{o} \mbox{ is bounded in } L^{1}\cap L^{\gamma}(\Omega), \qquad \rho^{n}_{o} \geq 0\mbox{ a.e. in }\Omega \\
\displaystyle\rho^{n}_{o}|u^{n}_{o}|^{2}=|m^{n}_{o}|^{2}/\rho^{n}_{o} \mbox{ is bounded in }  L^{1}(\Omega)\\ 
\displaystyle  \sqrt{\rho^{n}_{o}} \na \vphi(\rho^{n}_{o})=\na h(\rho^{n}_{o})/\sqrt{\rho^{n}_{o}} \mbox{ is bounded in } {L^{2}(\Omega)},\\
\displaystyle  \int \rho^{n}_{o} \frac{ |u^{n}_{o}|^{2+\delta}}{2} \, dx < C  \quad \mbox{ for some small $\delta$.}
\end{array}
\end{equation}
With those assumptions, and using inequalities (\ref{eq:entropy1}) and (\ref{eq:entropy2}), we deduce the following estimates, which we shall use
throughout  the proof of Theorem \ref{thm:stability}:
\begin{equation} \label{eq:bounds1n}
\begin{array}{l}
\ds || \sqrt \rho_{n} u_{n} ||_{L^{\infty}(0,T);L^{2}(\Omega))} \leq C \\
\ds || \rho_{n}||_{L^{\infty}(0,T;L^{1}\cap L^{\gamma}(\Omega))} \leq C \\
\ds ||\sqrt{h(\rho_{n})} \na u_{n}||_{L^{2}(0,T;L^{2}(\Omega))} \leq C
\end{array}
\end{equation}
and
\begin{equation}\label{eq:bounds2n}
\begin{array}{l}
\ds || h'(\rho_{n}) \na \sqrt \rho_{n}||_{L^{\infty}(0,T;L^{2}(\Omega))}\leq C\\
\ds ||\sqrt{h'(\rho_{n})\rho_{n}^{\gamma-2}}\na \rho_{n}||_{L^{2}(0,T;L^{2}(\Omega))}  \leq C
\end{array}
\end{equation}

\vspace{20pt}


In view of hypothesis on the viscosity coefficient (\ref{h}), the
bounds (\ref{eq:bounds1n}) and (\ref{eq:bounds2n}) yields:
\begin{equation}\label{eq:bounds3}
\begin{array}{l}
\ds ||\sqrt{\rho_{n}} \na u_{n}||_{L^{2}(0,T;L^{2}(\Omega))} \leq C\\
\ds ||  \na \sqrt \rho_{n}||_{L^{\infty}(0,T;L^{2}(\Omega))}\leq C\\
\ds || \na \rho_{n}^{\gamma/2}||_{L^{2}(0,T;L^{2}(\Omega))}  \leq C\\
\end{array}
\end{equation}
\vspace{20pt}

The proof of Theorem \ref{thm:stability} will be divided in 6 steps.
In the first two steps, we show the convergence of the density and
the pressure (note that  the convergence of the pressure is
straighforward here). The key argument of the proof is presented in
the third step: We prove that $\sqrt {\rho_{n}} u_{n}$ is bounded in
a space better than $L^{\infty}(0,T;L^{2}(\Omega))$. In turn, this
will give the convergence of the momentum (step 4) and finally the
strong convergence of $\sqrt{\rho_{n}} u_{n}$ in
$L^{2}_{loc}((0,T)\times \Omega)$ (step 5). The last step adresses
the convergence of the diffusion terms; It is mainly technical and
of minor interest. \vspace{10pt}

 \noindent{\bf Step 1: Convergence
of $\sqrt{\rho_{n}}$.}
\begin{lemma}\label{lem:sqrt}
If $h$ satisfies (\ref{h}),  then
\begin{eqnarray*}
&& \sqrt{\rho_{n}} \mbox{ is bounded in }   {L^{\infty}(0,T;H^{1}(\Omega))}  \\
&& \pa_{t} \sqrt {\rho_{n} }\mbox{ is bounded  in } L^{2}(0,T;H^{-1}(\Omega)) .
\end{eqnarray*}
As a consequence,  up to a subsequence, $\sqrt {\rho_{n}}$ converges
almost everywhere and  strongly  in
$L^{2}(0,T;L^{2}_{loc}(\Omega))$. We write
$$ \sqrt{\rho_{n}} \longrightarrow \sqrt \rho\qquad  \mbox { a.e and $L^{2}_{loc}((0,T)\times\Omega)$ strong.}$$
Moreover, $\rho_{n}$ converges to $\rho$ in $C^{0}(0,T;L^{3/2}_{loc}(\Omega))$.
\end{lemma}
{\it Proof.} The second estimate in (\ref{eq:bounds3}), together
with the conservation of mass
$||\rho_{n}(t)||_{L^{1}(\Omega)}=||\rho_{n,o}||_{L^{1}(\Omega)}$
gives the $L^{\infty}(0,T;H^{1}(\Omega))$ bound. Next, we notice
that
\begin{eqnarray*}
 \pa_{t}\sqrt{\rho_{n}} & =& - \frac{1}{2}\sqrt{\rho_{n}} \div u_{n} - u_{n}\cdot \na \sqrt{\rho_{n} } \\
& =&  \frac{1}{2}\sqrt{\rho_{n}} \div u_{n} - \div (u_{n} \sqrt{\rho_{n} })
\end{eqnarray*}
which yields the second estimate and, thanks to Aubin's Lemma, gives the strong convergence in $L^{2}_{loc}((0,T)\times \Omega)$.
\vspace{5pt}

Sobolev imbedding insures  that  $\sqrt{\rho_{n}}$ is bounded in
$L^{\infty}(0,T;L^{q}(\Omega))$ for $q\in[2,+\infty[$ if $N=2$ and
$q\in[2,6]$ if $N=3$. In either cases we deduce that $\rho_{n}$ is
bounded in $L^{\infty}(0,T;L^{3}(\Omega))$, and therefore
$$\rho_{n}u_{n} =\sqrt{\rho_{n}}\sqrt{\rho_{n}} u_{n}\mbox{ is bounded in } L^{\infty}(0,T;L^{3/2}(\Omega)).$$
The continuity equation thus yields $\pa_{t}\rho_{n}$  bounded in
$L^{\infty}(0,T,W^{-1,3/2}(\Omega))$. Moreover, since $\na\rho_{n} =
2\sqrt{\rho_{n}} \na \sqrt{\rho_{n}}$, we also have that $\na
\rho_{n} $ is bounded in  $L^{\infty}(0,T;L^{3/2}(\Omega))$, hence
the compactness of $\rho_{n} $  in $C([0,T];L^{3/2}_{loc}(\Omega))$.
\vspace{20pt}

\noindent{\bf Step 2: Convergence of the pressure}
\begin{lemma} \label{lem:pressure}
The pressure $\rho_{n}^{\gamma}$ is bounded in $L^{5/3}((0,T)\times
\Omega)$ when $N=3$ and $L^{r}((0,T)\times \Omega)$ for all $r\in
[1,2[$ when $N=2$. In particular, $\rho_{n}^{\gamma}$ converges to
$\rho^{\gamma}$ strongly in $L^{1}_{loc}((0,T)\times\Omega)$.
\end{lemma}
{\it Proof.}
Inequalities (\ref{eq:bounds3}) and (\ref{eq:bounds1n}) yield $\rho_{n}^{\gamma/2} \in L^{2} (0,T;H^{1}(\Omega))$.

When $N=2$, we deduce $ \rho_{n}^{\gamma/2} \in
L^{2}(0,T;L^{q}(\Omega))$ for all $q\in[2,\infty[$. So $ \rho
_{n}^{\gamma}$  is bounded in $L^{1}(0,T;L^{p}(\Omega))\cap
L^{\infty}(0,T;L^{1}(\Omega))$ for all $p\in[1,\infty[$, hence $
\rho _{n}^{\gamma}$ is bounded in $L^{r}((0,T)\times \Omega)$ for
all $r\in[1,2[$.

When $N=3$, we only get $ \rho_{n}^{\gamma/2} \in L^{2}(0,T;L^{6}(\Omega))$, or
$$\rho_{n}^{\gamma} \in L^{1}(0,T;L^{3}(\Omega)).$$
Since $ \rho_{n}^{\gamma}$ is bounded in $L^{\infty}(0,T;L^{1}(\Omega))$,
  H\"older inequality gives
$$||\rho_{n}^{\gamma}||_{L^{5/3}((0,T)\times\Omega)} \leq || \rho_{n}^{\gamma}||^{2/5}_{L^{\infty}(0,T;L^{1}(\Omega))}  \, ||\rho_{n}^{\gamma}||^{3/5}_{L^{1}(0,T;L^{3}(\Omega))} \leq C.$$
hence  $\rho_{n}^{\gamma}$ is bounded in $L^{5/3}((0,T)\times\Omega)$.

Since we already know that $\rho_{n}^{\gamma}$    converges almost everywhere to $\rho^{\gamma}$, those bounds yield the strong convergence  of $\rho_{n}^{\gamma}$ in $L^{1}_{loc}((0,T)\times\Omega)$.

\vspace{20pt}

\noindent{\bf Step 3: Bounds for $\sqrt{\rho_{n}}u_n$}

\begin{lemma}\label{lem:moment}
If $\gamma <3$, or if $N=3$, $\gamma\geq 3$ and (\ref{hg}) holds,
then
$$ \sqrt{\rho_{n}} u_{n} \mbox{ is bounded in } L^{\infty}(0,T;L^{2+2\alpha}(\Omega))$$
for some small $\alpha>0$
\end{lemma}

This Lemma  is really the corner stone of the stability result.
As a matter of fact, at this point, the main difficulty is to prove the strong convergence of  $\sqrt{\rho_{n}}u_n$ in $L^{1}(0,T;L_{loc}^{2}(\Omega))$.
A first consequence of Lemma \ref{lem:moment}, is that it will be enough to prove the convergence almost everywhere.
However, since we are only able to prove the convergence of the momentum $\rho_{n}u_{n}$
(see Step 4, which makes use of Lemma \ref{lem:moment} as well), we need to control
$\sqrt{\rho_{n}}u_n$ on the vacuum set $\{\rho(t,x)=0\}$ (and prove that it converges to zero almost everywhere).
And this fact also  will  be a consequence of Lemma \ref{lem:moment} (see Step 5).
\vspace{10pt}

{\it Proof.} The proof of Lemma \ref{lem:moment} relies on Lemma \ref{lem:estimate2}:
for $\delta$ small enough ($\delta\in(0,\nu/4)$), we have:
\begin{eqnarray}
&&\frac{d}{dt} \int \rho \frac{|u|^{2+\delta}}{2+\delta}\,dx + \frac{\nu}{2} \int h(\rho)|u|^{\delta} |\na u|^{2} \, dx \nonumber\\
&&\qquad\qquad  \leq
\left( \int \left(\frac{\rho^{2\gamma-\delta/2}}{h(\rho)}\right)^{2/(2-\delta)}\, dx\right)^{(2-\delta)/2} \left(\int \rho |u|^{2}\,dx \right)^{\delta/2}. \label{eq:moment'}
\end{eqnarray}
Using (\ref{eq:bounds1n}), we deduce:
\begin{equation*}
\frac{d}{dt} \int \rho \frac{|u|^{2+\delta}}{2+\delta}\,dx  \leq C
\left( \int
\left(\frac{\rho^{2\gamma-\delta/2}}{h(\rho)}\right)^{2/(2-\delta)}\,
dx\right)^{(2-\delta)/2} .
\end{equation*}
Condition (\ref{h}) yields $h(\rho)\geq \nu \rho$
and so
\begin{equation*}
\frac{d}{dt} \int \rho \frac{|u|^{2+\delta}}{2+\delta}\,dx  \leq C
\left( \int \left( \rho^{2\gamma-1-\delta/2}\right)^{2/(2-\delta)}\,
dx\right)^{(2-\delta)/2} .
\end{equation*}
Using Lemma \ref{lem:pressure}, we readily check that
the right hand side is bounded $L^{1}$ in time (for small $\delta$),
without any condition when $N=2$, and when $N=3$ under the condition that
$$ 2\gamma-1 < \frac{5}{3}\gamma,$$
which gives rise to the restriction  $\gamma < 3$.
In either cases, we deduce
\begin{equation*}
\frac{d}{dt} \int \rho \frac{|u|^{2+\delta}}{2+\delta}\,dx  \leq C .
\end{equation*}
and (\ref{initial22}) gives the lemma.
When $N=3$ and  $\gamma\geq3$ we need the extra hypothesis (\ref{hg})  to achieve the same result.
\vspace{10pt}

Finally, for $\alpha < \delta/2$, we have:
$$ \int (\rho |u|^{2} )^{1+\alpha}\,dx \leq  \left(\int \rho |u|^{2+\delta}\, dx \right)^{\frac{2+2\alpha}{2+\delta}} \left(\int \rho^{q(1+\alpha-(2+2\alpha)/
(2+\delta))} \, dx \right)^{\frac{1}{q}} $$ with
$q=(1-(2+2\alpha)/(2+\delta))^{-1}$, so that the exponent of $\rho$
goes to $1$ when $\alpha$ goes to zero. In particular, it is less
that $3$ for $\alpha$ small enough, and since
$\rho_{n}$ is bounded in $L^{\infty}(0,T;L^{3}(\Omega))$, we deduce
 Lemma \ref{lem:moment}.
\vspace{20pt}

\noindent{\bf Step 4: Convergence of the momentum}
\begin{lemma} \label{lem:ru}
Up to a subsequence, the momentum $m_{n} = \rho_{n}u_{n}$ converges strongly in $L^{2}(0,T;L^{1+\eps}_{loc}(\Omega))$  (for some positive $\eps$) and almost everywhere to some $m(x,t)$.
\end{lemma}
Note that we can already define $u(x,t) = m(x,t)/\rho(x,t)$ outside
the vacuum set $\{\rho(x,t) = 0\}$, but we do not know yet whether
$m(x,t)$ is zero on the vacuum set.
\\
{\it Proof.} We have
$$ \rho_{n}u_{n} = \sqrt{\rho_{n}} \sqrt{\rho_{n}} u_{n},$$
where $\sqrt{\rho_{n}}$ is bounded in
$L^{\infty}(0,T;L^{q}(\Omega))$ for $q\in[2,+\infty[$ if $N=2$ and $q\in[2,6]$ if $N=3$;
Since $\sqrt{\rho_{n}} u_{n} $ is bounded in $L^{\infty}(0,T;L^{2}(\Omega))$, we deduce that
 $$\rho_{n}u_{n} \mbox{  is bounded in }  L^{\infty}(0,T,L^{q}(\Omega)) \mbox{ for all $q\in[1,3/2]$}.$$
\vspace{5pt}

Next, we have
\begin{eqnarray*}
\pa_{i}(\rho_{n} {u_{n}}_{j}) & = & \rho_{n} \pa_{i} {u_{n}}_{j} + {u_{n}}_{j}\pa_{i} \rho_{n} \\
& = & \sqrt{\rho_{n}} \sqrt{\rho_{n}} \pa_{i} {u_{n}}_{j} + 2\sqrt{\rho_{n}} {u_{n}}_{j}\pa_{i} \sqrt{\rho_{n}}.
\end{eqnarray*}
Using Lemma \ref{lem:moment} and (\ref{eq:bounds3}), it is readily seen that the second term is bounded in $L^{\infty}(0,T;L^{1+\eps}(\Omega))$ for some small $\eps>0$, while the first term is bounded in
$L^{2}(0,T,L^{q}(\Omega))$ for all $q\in[1,3/2]$.
Hence
$$\na(\rho_{n} u_{n}) \mbox{ is bounded in } L^{2}(0,T;L^{1+\eps}(\Omega)).$$
In particular, we have
$$\rho_{n} u_{n} \mbox{ bounded in } L^{2}(0,T;W^{1,1+\eps}(\Omega)).$$

\vspace{5pt}

It remains to show that for every compact set $K\subset\Omega$,we have
\begin{equation}\label{eq:mom}
\pa_{t}(\rho_{n} u_{n}) \mbox{ is bounded in } L^{5/3}(0,T;W^{-2,3/2}(K)).
\end{equation}
As a matter of fact,  we observe that  $W^{1,3}_{0}(K)\subset  L^{1+1/\eps}(K)$ for small $\eps$ (for $N=2$ or $3$), and therefore
$$ L^{1+\eps}(K) \subset W^{-1,3/2}(K) \subset W^{-2,3/2}(K),$$
so (\ref{eq:mom}) together with Aubin's Lemma, yields the
compactness of $\rho_{n}u_{n}$ in $L^{2}(0,T;L^{1+\eps}(K))$.
\vspace{10pt}

To prove (\ref{eq:mom}), we use the  momentum equation (\ref{eq:euler2}), first noticing from Lemma
\ref{lem:pressure} and Lemma \ref{lem:moment}
that
\begin{eqnarray*}
\div \left(\sqrt{\rho_{n}} u_{n} \otimes\sqrt{\rho_{n}}u_{n} \right) &\in& L^{\infty}(0,T;W^{-1,1+\eps}(K)) \\
\na \rho_{n}^{\gamma} &\in& L^{5/3} (0,T;W^{-1,1+\eps}(K)),
\end{eqnarray*}

So we only have to check that $\na (h(\rho_{n})\na u_{n})$ and $\na( g(\rho_{n})\div u_{n})$ are bounded in $L^{\infty}(0,T;W^{-2,3/2}(K))$.
To that purpose, we write
\begin{equation}\label{eq:hu}
h(\rho_{n}) \na u_{n} = \na (h(\rho_{n}) u_{n} ) - u_{n}\na h(\rho_{n}),
\end{equation}
(and similarly with $g(\rho_{n})$).
The second term in (\ref{eq:hu}) is
$$ u_{n}\na h(\rho_{n}) = \sqrt{\rho_{n}} u_{n} \frac{\na h(\rho_{n})}{\sqrt{\rho_{n}}} = 2 \sqrt{\rho_{n}} u_{n} h'(\rho_{n}) \na \sqrt{\rho_{n}}$$ which is bounded in $L^{\infty}(0,T;L^{1+\eps}(\Omega))$ thanks to  (\ref{eq:bounds2n}) and Lemma \ref{lem:moment}.
The first term in (\ref{eq:hu}) can be rewritten
$$ \na [h(\rho_{n}) u_{n}]  = \na\left[ \frac{h(\rho_{n})}{\sqrt{\rho_{n}}} \sqrt{\rho_{n}} u_{n}\right],$$
which is bounded in $L^{\infty}(0,T;W^{-1,3/2}(\Omega))$ thanks to the following lemma:
\begin{lemma}\label{lem:h}
For all compact set $K$,  $h(\rho_{n})/\sqrt{\rho_{n}}$ and $g(\rho_{n})/\sqrt{\rho_{n}}$ are bounded in $L^{\infty}(0,T;L^{6}(K))$.
\end{lemma}
The proof of this Lemma is a bit technical in full generality and
will be postponed to Appendix \ref{app:gh}. However, note that, in
the particular case $h(\rho) = \nu \rho$, we have
$h(\rho_{n})/\sqrt{\rho_{n}} = \sqrt \rho_{n}$ and  Lemma
\ref{lem:h} follows straightforwardly from Lemma \ref{lem:sqrt}.

\vspace{5pt}

We deduce that   $h(\rho_{n}) \na u_{n}$  and $g(\rho_{n})\div u_{n}$ are bounded in
$$L^{\infty}(0,T;W^{-1,3/2}(K)+L^{1+\eps}(K)),$$
and since $L^{1+\eps}(K)) \subset W^{-1,3/2}(K)$ we can  conclude that
$h(\rho_{n}) \na u_{n}$ and $g(\rho_{n})\div u_{n}$  are bounded in $L^{\infty}(0,T;W^{-1,3/2}(K))$, which conclude the proof of Lemma \ref{lem:ru}.
\vspace{20pt}


\noindent{\bf Step 5: Convergence of $\sqrt{\rho_{n}} u_{n}$}

\begin{lemma}
The quantity $\sqrt{\rho_{n}} u_{n}$
converges strongly in $L^{1}$ and $L^{2}_{loc}((0,T)\times \Omega)$ to $m/\sqrt\rho $ (defined to be zero when $m=0$).

In particular,  we have $m(x,t)=0$ a.e. on $\{\rho(x,t)=0\}$ and there exists a function $u(x,t)$ such that
$$ m(x,t) = \rho(x,t) u(x,t)$$
(note that $u$ is not uniquely defined on the vacuum set $\{\rho(x,t) = 0\}$).
\end{lemma}
{\it Proof.}
First of all, since $m_{n}/\sqrt{\rho_{n}}$ is bounded in $L^{\infty}(0,T;L^{2}(\Omega))$, Fatou's lemma yields
$$ \int \liminf \frac{m_{n}^{2}}{\rho_{n}} \, dx < \infty .$$
In particular, we have $m(x,t)=0$ a.e. in $\{\rho(x,t)=0\}$, and
 if we define $m^{2}/\rho$ to be $0$ when $m=0$, we have
$$ \int  \frac{m^{2}}{\rho} \, dx < \infty .$$

Lemma \ref{lem:moment} implies that $\sqrt{\rho_{n}} |u_{n}|$ is bounded in $L^{\infty}(0,T;L^{2+2\alpha}(\Omega))$ for a small $\alpha$.
It is thus enough to prove the convergence almost everywhere, or in $L^{1}_{loc}((0,T)\times \Omega)$, to prove the strong convergence in $L^{2}_{loc}$.

\vspace{5pt}

First of all, we note that in  $\{\rho(x,t)\neq 0\}$, $\sqrt{\rho_{n}} u_{n}$ converges almost everywhere to $m/\sqrt \rho $.
So, if we denote the vacuum set by
$$ A= \{\rho(x,t)=0\},$$
we deduce
\begin{equation}\label{eq:vacuum}
 \sqrt{\rho_{n}} u_{n} 1_{\C A} \longrightarrow \frac{m}{\sqrt \rho}1_{\C A}\quad \mbox{ a.e.}.
 \end{equation}

To controle $\sqrt{\rho_{n}}u_{n}$ on the vacuum set, we introduce the set
$$B^{n}_{M} = \{\rho_{n} ^{1/(2+\delta)} |u_{n}| \geq M\},$$
for $M>0$.
We then cut  the $L^{1}$ norm as follows:
$$\int |\sqrt{\rho_{n}} u_{n} - \frac{m}{\sqrt\rho}|\, dx\, dt =  \int_{\C B^{n}_{M}\setminus A} \cdots + \int _{\C B_{M}^{n}\cap A} \cdots + \int_{B_{M}^{n}} \cdots
$$

The $L^{\infty}(0,T;L^{2}(\Omega))$ bound and (\ref{eq:vacuum}) gives the convergence of the first integral:
\begin{equation}\label{eq:in2}
 \int 1_{\C B_{M}^{n}\setminus A} |\sqrt{\rho_{n}} u_{n} - \frac{m}{\sqrt\rho} |\, dx\, dt  \longrightarrow 0.
 \end{equation}

Moreover, Lemma \ref{lem:moment} and Tchebychev's inequality yields
$$ |B^{n}_{M} | \leq \frac{C}{M^{2}},$$
and so
\begin{equation}\label{eq:in1}
\int_{B_{M}^{n}} |\sqrt{\rho_{n}} u_{n} -m/\sqrt\rho |   \leq \sqrt{|B_{M}^{n}|} \left( \int \rho_{n}|u_{n}|^{2}+ |m|^{2}/\rho\,dx \right)^{1/2} \leq \frac{C}{M}.
\end{equation}

Finally, on $\C B_{M}^{n} \cap A$,  we have
$$ |\sqrt{\rho_{n}} u_{n} | \leq M \rho_{n}^{1/2 - 1/(2+\delta)}\longrightarrow 0 \quad \mbox{ a.e. } ,$$
since $\rho_{n} \rightarrow 0$ a.e. and $1/2 - 1/(2+\delta)>0$.
So $1_{\C B_{M}^{n}\cap A} |\sqrt{\rho_{n}} u_{n}|$ converges almost everywhere to $0$.
In particular, the $L^{\infty}(0,T;L^{2}(\Omega))$ bound  yields
$$ \int 1_{B_{M}^{c}\cap A} |\sqrt{\rho_{n}} u_{n}| \, dx\, dt  \longrightarrow 0.$$
Since we defined $m/\sqrt\rho$ to be $0$ on $A$, we also have
$$ 1_{\C B_{M}^{n}\cap A}\, \frac{m}{\sqrt\rho}(x,t) = 0  \mbox{ a.e. } \forall n $$
hence
\begin{equation}\label{eq:in3}
\int 1_{\C B_{M}^{n}\cap A} |\sqrt{\rho_{n}} u_{n}-\frac{m}{\sqrt\rho} | \, dx\, dt  \longrightarrow 0.
\end{equation}

Putting (\ref{eq:in2}), (\ref{eq:in1}) and (\ref{eq:in3}) together,
we deduce
$$ \limsup_{n\rightarrow \infty} \int |\sqrt{\rho_{n}} u_{n} - \frac{m}{\sqrt\rho }|\, dx\, dt  \leq \frac{C}{M}$$
for all $M>0$, and so $\sqrt{\rho_{n}} u_{n} $ converges to
$\frac{m}{\sqrt\rho}$ in $L^{1}((0,T)\times\Omega)$ strong. The
lemma follows.
\vspace{50pt}

\noindent {\bf Step 6: Convergence of the diffusion terms}
\begin{lemma}
We have
$$ h(\rho_{n}) \na u_{n} \longrightarrow h(\rho) \na u \mbox{ in } \mathcal D'$$
and
$$ g(\rho_{n}) \div u_{n} \longrightarrow g(\rho) \div u \mbox{ in } \mathcal D'$$

\end{lemma}
{\it Proof.}
Let $\phi$ be a test function, then
\begin{eqnarray*}
 &&\int h(\rho_{n}) \na u_{n} \phi\, dx\, dt   =   - \int h(\rho_{n})  u_{n} \na \phi\, dx\, dt +\int  u_{n} \na h(\rho_{n})\phi\, dx\, dt \\
 && \qquad\qquad =   - \int \frac{h(\rho_{n})}{\sqrt \rho_{n}}  \sqrt{\rho_{n}} u_{n} \na \phi\, dx\, dt +\int  \sqrt{\rho_{n}} u_{n}  \frac{h'(\rho_{n})}{\sqrt{\rho_{n}}}\na \rho_{n} \phi\, dx\, dt
\end{eqnarray*}
Thanks to Lemma \ref{lem:h}, we know that $\frac{h(\rho_{n})}{\sqrt
\rho_{n}} \mbox{ is bounded in } L^{\infty
}(0,T;L^{6}_{loc}(\Omega))$. Moreover, since $h(\rho_{n})/\sqrt
\rho_{n}\leq \nu \sqrt{\rho_{n}}$, this term converges almost
everywhere to $h(\rho)/\sqrt \rho$ (defined to be zero on the vacuum
set). Therefore, it converges strongly in
$L^{2}_{loc}((0,T)\times\Omega)$; This is enough to prove  the
convergence of the first term.

Next, we note that
$$\frac{h'(\rho_{n})}{\sqrt{\rho_{n}}}\na \rho_{n}  = \na \psi (\rho_{n})$$
with $\psi'(\rho) = h'(\rho)/\sqrt\rho = \sqrt \rho\vphi'(\rho) $. Since
$$ \int |\na \psi(\rho) |^{2}\,dx = \int \rho |\na \vphi(\rho)|^{2}\,dx,$$
we have that $\na \psi(\rho_{n})$ is bounded in $L^{\infty}(0,T,
L^{2}(\Omega))$. Moreover, (\ref{eq:h2/3}) yields
$$ h'(\rho) \leq C \rho^{-1/2+\nu/3} \mbox{ when $\rho\leq 1$ }$$
and so
$$ \psi(\rho) \leq   C \rho^{\nu/3} \mbox{ when $\rho\leq 1$.}$$
Therefore, an argument similar to the proof of Lemma \ref{lem:h} shows that $\psi(\rho_{n})$ is bounded in $L^{\infty}(0,T;L^{6}_{loc}(\Omega))$. Since it converges almost everywhere ($\psi$ is a continuous function), it converges strongly in $L^{2}_{loc}((0,T)\times\Omega)$. It follows that
$$ \na \psi(\rho_{n}) \rightharpoonup \na \psi(\rho)\qquad L^{2}_{loc}((0,T)\times\Omega) \mbox{-weak }.
$$
A similar argument holds for $g(\rho_{n})\div u_{n}$ using the fact that $|g(\rho)|\leq C h(\rho)$ and $|g'(\rho)| \leq C h'(\rho).$

\section{Proof of Lemma \ref{lem:estimate}}\label{sec:estimate}
We conclude this paper by giving the proof of the estimate (\ref{eq:entropy2}).
To that purpose, we have to evaluate
$$ \frac{d}{dt} \int \left[ \frac 1 2 \rho |u|^{2} + \rho u\cdot \na \vphi(\rho) +\frac 1 2 \rho |\na\vphi(\rho)|^{2}\right]\, dx +\frac{d}{dt} \int \frac{1}{\gamma-1} \rho^{\gamma}\, dx.$$

\vspace{10pt}

{\bf Step 1:} First of all, we recall the usual entropy equality:
$$ \frac{d}{dt}\int \left[ \frac{1}{2}\rho |u|^{2}+\frac{1}{\gamma-1} \rho^{\gamma}\right] \, dx = -\int h(\rho) |\na u|^{2}\, dx  - \int g(\rho) |\div u|^{2} \, dx$$

\vspace{10pt}

{\bf Step 2:} Next,  (\ref{eq:euler1}) gives
\begin{eqnarray}
&&\int \rho \pa t \frac{|\na \vphi(\rho)|^{2}}{2}\, dx -\int \frac{|\na \vphi(\rho)|^{2}}{2} \div \rho u \, dx \nonumber\\
&&\qquad =  -\int \rho \na u : \na \vphi(\rho)\otimes \na
\vphi(\rho)\, dx
 + \int \rho^{2}\vphi'(\rho) \Delta \vphi(\rho) \div u \, dx\nonumber \\
&&\qquad\qquad + \int \rho [\na \vphi(\rho)]^{2} \div u\, dx
\nonumber
\end{eqnarray}
and so
\begin{eqnarray}
\frac{d}{dt}\int \rho \frac{|\na \vphi(\rho)|^{2}}{2}\, dx & = & -\int \rho \na u : \na \vphi(\rho)\otimes \na \vphi(\rho)\, dx \nonumber \\
&& + \int \rho^{2}\vphi'(\rho) \Delta \vphi(\rho) \div u \, dx\nonumber \\
&& + \int \rho [\na \vphi(\rho)]^{2} \div u \,dx \label{step2}
\end{eqnarray}

\vspace{10pt}

{\bf Step 3:} It remains  to evaluate the derivative of the cross-product:
\begin{eqnarray}
\frac{d}{dt}\int \rho u \cdot\na \vphi(\rho)\, dx &  = & \int \na \vphi(\rho) \cdot \pa_{t} (\rho u) \, dx+ \int \rho u \cdot \pa_{t} \na \vphi(\rho)\, dx \nonumber\\
& = & \int \na \vphi(\rho) \cdot \pa_{t} (\rho u) \, dx- \int \div(\rho u) \vphi'(\rho) \pa_{t} \rho\, dx\nonumber \\
& = & \int \na \vphi(\rho) \cdot \pa_{t} (\rho u) \, dx+ \int (\div(\rho u))^{2} \vphi'(\rho) \, dx.\label{step3}
\end{eqnarray}
Multiplying (\ref{eq:euler2}) by $\na \vphi(\rho)$, we get:
\begin{eqnarray*}
&&\int \na \vphi(\rho) \cdot \pa_{t} (\rho u) \, dx\\
   &&\qquad=  - \int (h(\rho)+g(\rho))\Delta\vphi(\rho) \div u\, dx+ \int \na u : \na \vphi(\rho)\otimes \na h(\rho) \, dx\\
&&\qquad\qquad -\int  \na \vphi(\rho)\cdot \na h(\rho) \div u \, dx - \int \na \vphi(\rho)  \cdot \na \rho^{\gamma} \, dx \\
&& \qquad\qquad- \int \na \vphi(\rho) \div(\rho u \otimes u)\,dx,
\end{eqnarray*}
where we used the fact that
$$
\int \na (g(\rho) \div u) \cdot \na \vphi(\rho) \, dx = -\int g(\rho) \Delta \vphi(\rho) \div u \, dx$$
and
\begin{eqnarray*}
&&\int \div(h(\rho)\na u)\cdot \na \vphi(\rho)\, dx\\
 &&\qquad=  \int \pa_{j}(h(\rho)\pa_{j}u_{i}))\pa_{i}\vphi(\rho) \,dx\\
&&\qquad =   \int \pa_{i}(h(\rho)\pa_{j}u_{i}))\pa_{j}\vphi(\rho) \,dx\\
&&\qquad =  \int\pa_{i} h(\rho) \pa_{j}u_{i} \pa_{j}\vphi(\rho)\,dx  -\int \pa_i u_{i} \pa_{j} h(\rho)\pa_{j}\vphi(\rho) \,dx\\
& &\qquad\qquad -\int \pa_i u_{i} h(\rho) \pa_{jj} \vphi (\rho) \,dx\\
&&\qquad =  \int \na u : \na h(\rho)\otimes \na \vphi(\rho)\, dx
- \int \na h(\rho)\cdot \na \vphi(\rho)\div u\, dx\\
& &\qquad\qquad - \int  h(\rho) \Delta \vphi(\rho) \div u \, dx
\end{eqnarray*}

\vspace{10pt}

{\bf Step 4:} When $\vphi$, $h$ and $g$ satisfies (\ref{gh}) and
(\ref{phi}), then (\ref{step2}) and (\ref{step3}) yields
\begin{eqnarray*}
&&  \frac{d}{dt}\left\{ \int\rho u \cdot\na \vphi(\rho) + \rho \frac{|\na \vphi(\rho)|^{2}}{2}\, dx \right\} +
\int \na \vphi(\rho) \cdot \nabla p\, dx \\
&& \qquad\qquad =- \int \na \vphi(\rho) \div(\rho u \otimes u)\, dx +  \int \vphi'(\rho) (\div(\rho u))^{2} \, dx.
\end{eqnarray*}
Finally, we have
\begin{eqnarray*}
&& - \int \na \vphi(\rho) \div(\rho u \otimes u) \, dx+  \int \vphi'(\rho) (\div(\rho u))^{2} \, dx \\
&& \qquad = \int - \vphi'(\rho) u\cdot \na \rho \, \div (\rho u) -\vphi'(\rho) \na \rho (\rho u \cdot\na u) + \vphi'(\rho) (\div \rho u)^{2} \, dx\\
&& \qquad = \int \rho \vphi'(\rho)  \, \div  u \, \div (\rho u) - \rho\vphi'(\rho) \na \rho (u\cdot\na u) \, dx\\
&& \qquad = \int \rho^{2}\vphi'(\rho)  (\div  u )^{2} + \rho \vphi'(\rho) u\cdot \na \rho \, \div u - \rho \vphi'(\rho) \na \rho (u \cdot\na u) \, dx
\end{eqnarray*}
so using (\ref{phi}) and (\ref{gh}), we get
\begin{eqnarray*}
&& - \int \na \vphi(\rho) \, \div(\rho u \otimes u)\, dx +  \int \vphi'(\rho) (\div(\rho u))^{2} \, dx\\
&& \quad = \int \rho h'(\rho)  (\div  u )^{2} + \na(h(\rho)) \cdot u \, \div u - \na(h( \rho)) (u \cdot\na u) \, dx\\
&& \quad= \int \rho h'(\rho) (\div u)^{2} - h(\rho) (\div u)^{2}- h(\rho) u\cdot \na\div u\,dx\\
&&\qquad\qquad + \int h(\rho) \pa_{i } u _{j} \pa_{j}u_{i} + h(\rho) u \cdot \na \div u \, dx\\
&& \quad = \int (\rho h' - h )( \div u)^{2} + h(\rho) \pa_{i } u _{j} \pa_{j} u_{i}\, dx \\
&& \quad = \int g(\rho) \, (\div u)^{2}\, dx + \int h(\rho) \pa_{j}u_{i} \pa_{i}u_{j}\, dx\\
\end{eqnarray*}
which yields
\begin{eqnarray*}
&&  \frac{d}{dt}\left\{ \int\rho u \cdot\na \vphi(\rho) + \rho \frac{|\na \vphi(\rho)|^{2}}{2}\, dx \right\} +
\int \na \vphi(\rho) \cdot \nabla \rho^{\gamma}\, dx \\
&& \qquad\qquad \leq \int g(\rho) (\div u)^{2}\, dx + \int h(\rho) |\na u|^{2}\, dx,
\end{eqnarray*}
and the proof is complete.

\appendix

\section{Proof of Lemma \ref{lem:h}} \label{app:gh}
We shall only prove the result for $h(\rho_{n})/\sqrt{\rho_{n}}$.
Using the fact that
$$|g(\rho)|\leq C h(\rho), \quad \mbox{and } |g'(\rho)|\leq C h'(\rho) \quad \mbox{ for all } \rho,$$
a similar proof follows for $g(\rho_{n})/\sqrt{\rho_{n}}$

Note that In view of (\ref{eq:h2/3}), we have
$$ \frac{h(\rho)}{\sqrt{\rho}} \leq C\rho^{\nu} \qquad \mbox{ if } \rho\leq 1 ,$$
so we only need to  control  $\frac{h(\rho_{n})}{\sqrt{\rho_{n}}}$ for large $\rho_{n}$.
This will be achieved differently depending on the dimension.
\vspace{5pt}

When $N=2$, the fact that $\sqrt{\rho_{n}}$ is bounded in
$L^{\infty} (0,T;H^{1}(\Omega))$ and Sobolev's inequalities implies
that $\rho_{n}$ is bounded in $L^{\infty} (0,T;L^{p}(\Omega))$ for
all $p\in[1,\infty[$. Moreover, in view of  (\ref{eq:h2/3}), we have
$$ \frac{h(\rho)}{\sqrt{\rho}} \leq
\left\{\begin{array}{ll}
C\rho^{1/\nu} & \mbox{ if } \rho\geq 1 \\
C\rho^{\nu} & \mbox{ if } \rho\leq 1
\end{array}\right. .$$
So there exists $q_{o}>1$ such that $h(\rho_{n})/\sqrt{\rho_{n}}$ is bounded in $L^{\infty}(0,T;L^{q}(\Omega))$ for all $q > q_{o}$. In particular, $h(\rho_{n})/\sqrt{\rho_{n}}$ is bounded in $L^{\infty}(0,T;L^{p}(K))$ for all $p\in[1,\infty[$ for any compact set $K$.
\vspace{5pt}

When $N=3 $, we note that
$$ \na\left(\frac{h(\rho)}{\sqrt \rho}\right) = 2h'(\rho)\na \sqrt \rho - \frac{h(\rho)}{2\rho^{3/2}} \na \rho,$$
and since conditions (\ref{gh}) and (\ref{g+h}) yields
$$  h'(\rho)\rho =g(\rho)+h(\rho)\geq  \frac{3 g (\rho) +h(\rho)}{3}  \geq  \frac{\nu}{3} h(\rho),$$
we have
$$ |\na\left(\frac{h(\rho)}{\sqrt \rho}\right)| \leq C |h'(\rho)\na \sqrt \rho|.$$
So inequality (\ref{eq:bounds2}) yields
\begin{equation} \label{eq:bound4}
||\na\left(\frac{h(\rho_{n})}{\sqrt{ \rho_{n}}}\right)||_{L^{\infty}(0,T;L^{2}(\Omega))} \leq C
\end{equation}
When $\Omega=\R^{3}$, Sobolev's inequalities implies  that
$h(\rho_{n})/\sqrt{\rho_{n}}$ is bounded in
$L^{\infty}(0,T;L^{6}(\Omega))$. When $\Omega$ is a subset of
$\R^{3}$, we note that (\ref{eq:h2/3}) gives
$$ \frac{h(\rho)}{\sqrt{\rho}} \leq
\left\{\begin{array}{ll}
C\rho^{1/6+3/\nu} & \mbox{ if } \rho\geq 1 \\
C\rho^{1/6+\nu/3} & \mbox{ if } \rho\leq 1
\end{array}\right. .$$
So there exists a constant $s\leq1$ such that
$$ \left( \left(\frac{h(\rho_{n})}{\sqrt{\rho_{n}}}\right)^{s}-1 \right)_{+} \in L^{\infty}(0,T;L^{2}(\Omega))$$
Moreover
\begin{eqnarray*}
 \left| \na \left(\frac{h(\rho_{n})}{\sqrt{\rho_{n}}}\right)^{s} 1_{h(\rho_{n})/\sqrt{\rho_{n}}\geq 1}  \right|& = &
\left | \left(\frac{h(\rho_{n})}{\sqrt{\rho_{n}}}\right)^{s-1}
\na \left(\frac{h(\rho_{n})}{\sqrt{\rho_{n}}}\right) 1_{h(\rho_{n})/\sqrt{\rho_{n}}\geq 1} \right| \\
& \leq &   \left|\na \left(\frac{h(\rho_{n})}{\sqrt{\rho_{n}}}\right)\right|\in L^{\infty}(0,T;L^{2}(\Omega)),
\end{eqnarray*}
using the fact that $s-1\leq0$.
It follows that $(h(\rho_{n})/\sqrt{\rho_{n}})^{s} 1_{\rho_{n}\geq 1}$ is bounded in $L^{\infty}(0,T;H^{1}(\Omega))$ which in turn gives
$$ \left(\frac{h(\rho_{n})}{\sqrt{\rho_{n}}}\right)^{s_{1}} 1_{\rho_{n}\geq 1} \in L^{\infty}(0,T;L^{2}(\Omega)),$$
for all $s_{1}\in(s,3s)$.
As long as $3s \leq 1$, we can repeat this argument with $3s$ instead of $s$.
Eventually, this will lead to
$$ \left(\frac{h(\rho_{n})}{\sqrt{\rho_{n}}}\right) 1_{\rho_{n}\geq 1} \in L^{\infty}(0,T;L^{2}(\Omega)),$$
which, together with (\ref{eq:bound4}) implies that $(h(\rho_{n})/\sqrt{\rho_{n}}) 1_{\rho_{n}\geq 1}$ is bounded in $L^{\infty}(0,T;L^{6}(\Omega))$.

\bibliography{biblio}

\end{document}